\documentclass{article}
\usepackage{hyperref,amsmath,amsthm,amssymb,setspace}
\usepackage[dvips]{graphicx}

\title{Variations of a Coin-Removal Problem}
\author{Kennan Shelton\\
   {\small Mathematics and Computer Science}\\[-.4em]{\small Rhodes College}\\[-.4em]{\small Memphis, TN 38112}\\[-.4em]{\small \tt{shelton@rhodes.edu}}
\and Michael Siler\\
{\small Box 1959}\\[-.4em]{\small Rhodes College}\\[-.4em]{\small
Memphis, TN 38112}\\[-.4em]{\small \tt{silwm@rhodes.edu}}
}
\date{October, 2004}

\newtheorem{thm}{Theorem}
\newtheorem*{thm*}{Theorem}
\newtheorem{cor}{Corollary}
\newtheorem{lem}{Lemma}
\newtheorem{prop}{Proposition}
\newtheorem{conj}{Conjecture}
\newtheorem*{conj*}{Conjecture}
\theoremstyle{definition}

\theoremstyle{remark}
\newtheorem{rem}{Remark}

\newcommand{\rh}[1]{\overset{\bullet}{#1}}
\newcommand{\flip}[1]{\overline{#1}}

\begin{document}

\maketitle

\begin{abstract}
Given a set of coins arranged in a line, we remove heads-up coins one at a time
and flip any adjacent coins after each removal.  The coin-removal problem is to
determine for which arrangements of coins it is possible to remove all of the
coins.  In this paper we consider a variation of the problem in which gaps
created by removing coins are eliminated by pushing the coins together. We
characterize the set of removable arrangements and show that this set forms a
regular language.  We use a finite automaton to find a recursive formula for
the number of removable arrangements of different lengths.
\end{abstract}

\section{Introduction}

The coin-removal problem is given in \cite[p. 50]{dangelowest}:
\begin{quote}
Suppose that $n$ coins are arranged in a row.  We remove heads-up coins, one by
one.  Each time we remove a coin we must flip the coins still present in the
(at most) two positions surrounding it.  For which arrangements of heads and
tails can we remove all the coins?
\end{quote}
We will call an arrangement {\em removable} if we can remove all the coins; the
coin-removal problem is to characterize the removable arrangements of coins. We
will also use the word ``sequence'' to refer to a linear arrangement of coins.

For example, if we use $1$ to denote a heads-up coin and $0$ to denote a
tails-up coin, then the sequence $01110$ is removable, as shown by the
following sequence of moves (the coin to be removed at a particular stage is
marked by a dot, an empty position by a circle):

\begin{equation*}
 0\rh{1}110 \rightarrow \rh{1}\circ010
        \rightarrow \circ\circ0\rh{1}0
        \rightarrow \circ\circ\rh{1}\circ1
        \rightarrow \circ\circ\circ\circ\rh{1}
        \rightarrow \circ\circ\circ\circ\circ
\end{equation*}
The solution of the coin-removal problem is simple to state
\cite[p.63]{dangelowest}: a sequence is removable if and only if the number of
heads-up coins is odd. Proving this statement was a homework assignment in a
recent ``Introduction to Proofs'' class taught at Rhodes College---the proof
found by the second author is presented in the next section.

It is natural to ask how the coin-removal problem might be modified or
generalized. The variation we consider in this paper is to use {\em linear
arrangements with no gaps}---any gap created by removing a coin is eliminated
by pushing the remaining coins together.  We give a  method for determining the
removability of a sequence in this variation.  Using a deterministic finite
automaton, we also show that the set of removable sequence in this variation
forms a regular language and find a recursive formula for the number of
removable sequences of different lengths. Finally, we discuss circular
arrangements of coins and suggest possible directions for further
investigations.

We will use the notation $1^n$ to denote a sequence of $n$ $1$'s (heads-up
coins). If $n = 0$ then the sequence is empty and any adjacent sequences can be
concatenated. For example, $101^001$ is the same sequence as $1001$.  Similar
notation will be used for sequences of tails-up coins.

\section{The Original Problem---Linear Arrangement With Gaps}

In the original formulation  of the coin-removal problem, the coins are
arranged in a line and removing coins may create gaps.  As we saw above, the
sequence $01110$ is removable. Other removable sequences include $1$, $111$ and
$11000111$. Clearly the sequence $11$ is not removable, nor are the sequences
$101$, $11011$, and $011011$.  Note that by allowing gaps, any remaining
non-adjacent sequences of coins are isolated and removing coins from one will
not affect the others.

\begin{prop}
\label{prop:linegaps} In the original coin-removal problem (linear arrangement
with gaps), an arrangement of coins is removable if and only if the number of
heads-up coins is odd.
\end{prop}
\begin{proof}
We will use strong mathematical induction on the length $k$ of the sequence of
coins.  Clearly the result holds for $k = 1$.

Assume $k\geq 1$ and that our claim holds for all sequences of lengths $1$ to
$k$.  Let $A = a_1 a_2\dots a_{k+1}$ be a coin sequence of length $k + 1$ (each
$a_n$ is $0$ or $1$).  Suppose first that $A$ has an even number of heads-up
coins (i.e., $a_n = 1$ for an even number of terms).  If there are zero
heads-up coins then $A$ is clearly not removable.  Suppose $A$ has at least two
heads-up coins.  Select any heads-up coin in $A$, say at position $n+1$, $0
\leq n \leq k$, so that $A = a_1\dots a_{n} 1 a_{n+2}\dots a_{k+1}$.  Then
either subsequence $a_1\dots a_n$ or $a_{n+2}\dots a_{k+1}$ contains an odd
number of heads-up coins.  Removing the selected coin flips any adjacent coins
so this subsequence will now contain an even number of heads-up coins.  Its
length is less than $k+1$ so by the inductive hypothesis, it is not removable.
Our choice of heads-up coin to remove was arbitrary and so our original
sequence is not removable.

Now suppose $A$ has an odd number of heads.  Then $A = 0^n 1 D$ where $n \geq
0$ and $D$ is a possibly empty sequence containing an even number of heads.
Removing the first heads-up coin on the left leaves us with one sequence or two
isolated sequences, each of length less than $k+1$ and containing an odd number
of heads-up coins.  By the inductive hypothesis, these two sequences are
removable, and so the original sequence is removable.

\end{proof}

\section{The Variation---Linear Arrangement With No Gaps}

Now suppose that the coins are again arranged a line but when a coin is
removed, any gap is eliminated by pushing the remaining coins together.  We
will say that a sequence $A$ is {\em reducible} to another sequence $A'$ if,
starting with the configuration $A$, we can remove heads-up coins and eliminate
gaps to obtain $A'$.

For example, consider the sequence of coins $11101$.  With gaps this sequence
is not removable, but without gaps it is reducible to $1$ and thus removable:
\begin{equation*}
1 1 \rh{1} 0 1 \rightarrow 1 0 \rh{1} 1 \rightarrow 1 \rh{1} 0 \rightarrow 0
\rh{1} \rightarrow \rh{1} \rightarrow \circ.
\end{equation*}
Also, the sequence $11111$ is removable with gaps but not in this variation.
Other examples of removable sequences include $1$, $01$, $111$, $1101$ and
$010111$. Examples of non-removable sequences include $11$, $010$, $0111$ and
$110110$.

Our approach will be to find an {\em invariant}, a property of removable
sequences that is unchanged under the removal of heads-up coins.  We will then
show that this property can distinguish between removable and non-removable
sequences.  It may not be immediately obvious why the invariant we use is
indeed a useful one, so we start with some motivation for its definition.

We first notice that sequences $1$, $111$, $1111$, $111111$ and $1111111$ are
removable while the sequences $11$, $11111$ and $11111111$ are not.  Based on
this (admittedly meager) evidence, we form the following conjecture:

\begin{conj}
The sequence $1^n$ is removable if and only if $n$ is congruent to $0$ or $1$
modulo $3$.
\end{conj}

Our conjecture suggests that we may be able to decide removability of a given
sequence $A$ by removing heads-up coins to reduce it to a sequence $1^n$ of
only heads.  If $n$ is congruent to $0$ or $1$ modulo $3$, and our conjecture
is true, then the sequence is removable.  Of course, to conclude that a
sequence is not removable would still require that we check every possible way
to reduce $A$ to a sequence $1^n$ and verify that in each case that $n$ is
congruent to $2$ modulo $3$.

With our conjecture in mind we next make the following observations:
\begin{enumerate}
\item A sequence of the form $1^m 00 1^n$, $m \geq 2$, $n\geq 0$, can be reduced to the sequence
$1^{m+n}$:
\begin{equation*}
1^m001^n = 1^{m-2}1\rh{1}001^n \rightarrow 1^{m-2}0\rh{1}01^n \rightarrow
1^{m-2}111^n = 1^{m+n}.
\end{equation*}
The same result holds when $m = 0$ or $m = 1$.

\item A sequence of the form $1^m01^nC$, $m\geq n\geq 1$, $C$ any sequence, can be reduced to the
sequence $1^{m-n}0C$:
\begin{align*}
1^m01^nC = 1^{m-2}1\rh{1}01^nC &\rightarrow 1^{m-2}0\rh{1}1^nC\\
            &\rightarrow 1^{m-2}101^{n-1}C = 1^{m-1}01^{n-1}C\\
            &\phantom{\rightarrow}\vdots\\
            &\rightarrow 1^{m-n+1}0\rh{1}C\\
            &\rightarrow 1^{m-n+1}\rh{1}C' \qquad (\text{$C'$ is $C$ with the first coin flipped})\\
            &\rightarrow 1^{m-n}0C.
\end{align*}

\item A sequence of the form $1^m0$, $m \geq 2$, can be reduced to the sequence
$1^{m-1}$:
\begin{equation*}
1^m0 = 1^{m-2}1\rh{1}0 \rightarrow 1^{m-2}0\rh{1} \rightarrow 1^{m-2}1 =
1^{m-1}.
\end{equation*}
The same result holds for $m = 1$ if we interpret $1^0$ as the empty sequence:
all coins have been removed.
\end{enumerate}

In general, two groups of $1$'s separated by an even group of $0$'s may be
reduced to a single group of $1$'s whose length is the {\em sum} of the lengths
of the original groups.  A sequence consisting of two groups of $1$'s separated
by an odd group of $0$'s may be reduced to a single group of $1$'s whose length
is the {\em difference} of the lengths of the original groups, followed by a
$0$.

For example, by applying the observations to the sequence $11110011111011001$,
we have that
\begin{alignat*}{2}
11110011111011001 &\rightarrow 111111111011001  & &\quad \text{(observation 1)}\\
            &\rightarrow 11111110001      & &\quad \text{(observation 2)}\\
            &\rightarrow 111111101          & &\quad \text{(observation 1, applied
            to $1111111001^001$)}\\
            &\rightarrow 1111110           & &\quad \text{(observation 2)}\\
            &\rightarrow 11111 = 1^5        & &\quad \text{(observation 3)}
\end{alignat*}
Since $5 \equiv 2 \pmod{3}$, we conjecture that $11110011111011001$ is not
removable.

The above observations together with the conjecture suggest that we may be able
to determine removability by simply adding or subtracting the numbers of $1$'s
in each group of $1$'s in the sequence (starting from the left) and examining
the result modulo $3$.  We will make this approach more precise and then show
that it does indeed provide a means of classifying sequences.

Given a coin sequence $A$ containing at least one $1$, we assume that $A$
starts and ends with (possibly empty) blocks of $1$'s (heads-up coins).  Then
\begin{equation*}
A = 1^{h_0} 0^{t_1} 1^{h_1} 0^{t_2}\cdots 1^{h_{n-1}} 0^{t_n} 1^{h_n}
\end{equation*}
for some choice of non-negative integer sequences $H(A) = (h_0,h_1,\dots,h_n)$
and $T(A) = (t_1,\dots t_n)$, $n\geq 1$.  Note that the choice of $H(A)$ and
$T(A)$ is not unique. For example, the sequence $A = 11011$ may be expressed as
$H(A) =(2,2)$ and $T(A) = (1)$ or $H(A) = (2,1,1)$ and $T(A) = (1,0)$.

Suppose that $H(A)$ and $T(A)$ have been chosen and let $P(A) =
(p_1,p_2,\dots,p_n)$ where $p_i = (t_1 + t_2 + \dots + t_n) \pmod{2}$.  The
sequence $P(A)$ keeps track of the modified parities of the groups of $0$'s: by
observation 2, if we have a group of $0$'s of odd parity then the parity of the
{\em next} group of $0$'s will change.

Finally, we define the quantity $S(A)$ by the following formula:
\begin{equation*}
S(A) = h_0 + \sum_{i=1}^n (-1)^{p_i}h_i - p_n.
\end{equation*}
We will call $S(A)$ the {\em parity sum} for $A$.  The parity sum adds or
subtracts the numbers of $1$'s in each group followed by a final subtraction of
$1$ if $A$ reduces to a sequence of $1$'s followed by a single $0$ (observation
3).

For example, for $A = 110001101110$, we can set $H(A) = (2,2,3,0)$ and $T(A) =
(3,1,1)$.  Then $P(A) = (1,0,1)$ and $S(A) = 2 - 2 + 3 - 0 - 1 = 2$.  Thus we
conjecture that $A$ is reducible to $11$ and so is not removable.

\begin{rem}
While we will omit its proof, it is not hard to show that the parity sum $S(A)$
is independent of the choice of sequences $H(A)$ and $T(A)$.  This is not true
for sequences of all $0$'s, hence the restriction of $S$ to only those with at
least one $1$.
\end{rem}

The parity sum will provide us with the invariant we need.  In particular, a
sequence $A$ containing at least one heads-up coin has the property that $S(A)$
is congruent to $2$ modulo $3$ if and only if the parity sum of any sequence
reducible from $A$ is also congruent to $2$ modulo $3$.  This property of the
parity sum will allow us to distinguish between removable and non-removable
sequences.

\begin{lem}
\label{lem:revSA} If $A$ is a coin sequence containing at least one heads-up
coin and $A^R$ is the reversal of $A$ then $S(A) \equiv 2 \pmod{3}$ if and only
if $S(A^R) \equiv 2 \pmod{3}$.
\end{lem}

\begin{proof}
Let $H(A) = (h_0,\dots,h_n)$ and $T(A)= (t_1,\dots,t_n)$ be fixed and set
$H(A^R) = (h_n,\dots,h_0)$ and $T(A) = (t_n,\dots,t_n)$. Then $P(A^R) =
(q_1,q_2,\dots,q_n)$ where $q_i = (t_n + t_{n-1} + \dots + t_{n-i+1})
\pmod{2}$.  Thus $p_n = q_n = (p_i + q_{n-i}) \pmod{2}$ and so
\begin{align*}
S(A) & = h_n + \sum_{i=1}^n (-1)^{q_i} h_{n-i} - q_n\\
     & = h_n + \sum_{i=0}^{n-1}(-1)^{q_{n-i}}h_i - q_n\\
     & = (-1)^{q_n}h_0 + \sum_{i=1}^{n-1}(-1)^{p_n - p_i}h_i + h_n - p_n\\
     & = (-1)^{p_n}h_0 + (-1)^{p_n}\sum_{i=1}^{n-1}(-1)^{p_i}h_i + h_n - p_n.
\end{align*}
If $p_n = 0$ then $S(A^R) = S(A)$.  If $p_n = 1$ then $S(A^R) = -h_0 -
\sum_{i=1}^n (-1)^{p_i}h_i - p_n = -S(A) - 2$.  Thus $S(A^R) + S(A) = -2 \equiv
1 \pmod{3}$.  In either case, $S(A) \equiv 2\pmod{3}$ if and only if $S(A^R)
\equiv 2 \pmod{3}$.
\end{proof}

\begin{lem}
\label{lem:SAequiv} Suppose $A$ is a coin sequence containing at least one
heads-up coin and $A'$ a sequence obtained from $A$ by removing a heads-up
coin.  Suppose further that $A'$ also contains at least one heads-up coin.
Then $S(A) \equiv 2 \pmod 3$ if and only if $S(A') \equiv 2 \pmod 3$.
\end{lem}

\begin{proof}
Without loss of generality, we may consider only the five cases listed below
for the original sequence $A$ and reduced sequence $A'$. In each case the coin
to be removed is marked by a dot; $C$ and $D$ are (possibly empty) sequences of
$0$'s and $1$.
\begin{enumerate}
    \item $A = \rh{1}0D \rightarrow A' = 1D$
    \item $A = \rh{1}1D \rightarrow A' = 0D$; $D$ contains a $1$
    \item $A = C0\rh{1}0D \rightarrow A' = C11D$
    \item $A = C1\rh{1}0D \rightarrow A' = C01D$
    \item $A = C1\rh{1}1D \rightarrow A' = C00D$; at least one of $C$ or $D$ contains
    a $1$.
\end{enumerate}
All other ways that we can remove a heads-up coin can be reduced to these cases
using the reversal of $A$ and Lemma~\ref{lem:revSA}.  For each case, we will
find a relationship between $S(A)$ and $S(A')$.

{\bf Case 1:} Assume that $H(A) = (h_0,h_1,\dots,h_n)$ with $h_0 = 1$. After
removing the indicated coin we can set $H(A') = H(A)$ and $T(A') = (t_1 -
1,t_2,\dots,t_n)$. Then $P(A') = (\flip{p_1},\flip{p_2},\dots,\flip{p_n})$
where $\flip{p} = 1 - p$. So
\begin{align*}
S(A') & = h_0 + \sum_{i=1}^n (-1)^{\flip{p_i}}h_i - \flip{p_n}\\
      & = 1 - \sum_{i=1}^n (-1)^{p_i}h_i + p_n - 1\\
      & = 1 - S(A).
\end{align*}
Thus $S(A') + S(A) = 1$ in Case 1.

{\bf Case 2:} Assume that $H(A) = (h_0,h_1,\dots,h_n)$ with $h_0 \geq 2$ and
$T(A) = (t_1,t_2,\dots,t_n)$. After removing the coin we can use the sequences
$H(A') = (0,h_0 - 2, h_1, \dots,h_n)$ and $T(A') = (1,t_1,t_2,\dots,t_n)$.  So
$P(A') = (1,\flip{p_1},\flip{p_2},\dots,\flip{p_n})$ and
\begin{align*}
S(A') & = 0 - (h_0 - 2) + \sum_{i=1}^{n}(-1)^{\flip{p_i}}h_i - \flip{p_n}\\
      & = 2 - h_0 - \sum_{i=1}^{n} (-1)^{p_i}h_i + p_n - 1\\
      & = 1 - S(A).
\end{align*}
So $S(A') + S(A) = 1$ in Case 2 as well.

{\bf Case 3:} Suppose the coin to be removed is in the $j^{\text{th}}$ group of
heads in $H(A)$, $1 \leq j < n$, so $h_j = 1$.  Then $H(A') =
(h_0,h_1,\dots,h_{j-1},h_j + 1,h_{j+1},\dots,h_n)$ and $T(A') =
(t_1,\dots,t_{j-1},t_j - 1,t_{j+1}-1,t_{j+2},\dots,t_n)$.  Thus $P(A') =
(p_1,\dots,\flip{p_j},p_{j+1},\dots,p_n)$ and so
\begin{align*}
S(A') & = h_0 + \sum_{i=1}^{j-1} (-1)^{p_i}h_i + (-1)^{\flip{p_j}}(h_j + 1) +
\sum_{i=j+1}^n (-1)^{p_i}h_i - p_n\\
      & = h_0 + \sum_{i=1}^{j-1} (-1)^{p_i}h_i + (-1)^{p_j}(-2) +
      \sum_{i=j+1}^n (-1)^{p_i}h_i - p_n\\
      & = h_0 + \sum_{i=1}^n (-1)^{p_i}h_i + (-1)^{p_j}(-3)\\
      & = S(A) - 3(-1)^{p_j}.
\end{align*}
Therefore $S(A') \equiv S(A) \pmod{3}$.

{\bf Case 4:} Suppose the coin to be removed is in the $j^{\text{th}}$ group of
heads in $H(A)$, $0 \leq j < n$, so $h_j \geq 2$.  Then we may take $H(A') =
(h_0,h_1,\dots,h_{j-1},h_j - 2,1,h_{j+1},\dots,h_n)$ and $T(A') =
(t_1,t_2,\dots,t_j,1,t_{j+1}-1,t_{j+2},\dots,t_n)$.  From this we have that
$P(A') = (p_1,\dots,p_j,\flip{p_j},p_{j+1},\dots,p_n)$ (if $j = 0$ then $P(A')
= (1,p_1 -1,p_2,\dots,p_n)$. Then
\begin{align*}
S(A') & = h_0 + \sum_{i=1}^{j-1}(-1)^{p_i} h_i + (-1)^{p_j}(h_j - 2) +
(-1)^{\flip{p_j}}\cdot 1 + \sum_{i=j+1}^n (-1)^{p_i} h_i - p_n\\
      & = S(A) - 3(-1)^{p_j}
\end{align*}
and again we have that $S(A') \equiv S(A) \pmod{3}$.

{\bf Case 5:} Suppose the coin to be removed is in the $j^{\text{th}}$ group of
heads in $H(A)$, $0 \leq j < n$, so $h_j \geq 3$.  Then take $H(A') =
(h_0,h_1,\dots,h_j - 3,0,h_{j+1},\dots,h_n)$ and $T(A') =
(t_1,\dots,t_j,2,t_{j+1},\dots,t_n)$, so $P(A') =
(p_1,\dots,p_j,p_j,p_{j+1},\dots,p_n)$.  Thus
\begin{align*}
S(A') & = h_0 + \sum_{i=1}^{j-1}(-1)^{p_i}h_i + (-1)^{p_j}(h_j - 3) +
(-1)^{p_j}\cdot 0 + \sum_{i=j+1}^n (-1)^{p_i} h_i - p_n\\
      & = S(A) - 3(-1)^{p_j}
\end{align*}
and so $S(A') \equiv S(A) \pmod{3}$.

Therefore, in each case we have that $S(A) \equiv 2 \pmod{3}$ if and only if
$S(A') \equiv 2 \pmod{3}$ and this property of $S(A)$ is invariant under the
removal of a heads-up coin.

\end{proof}

The following corollary is an immediate consequence of our previous
observations and Lemma~\ref{lem:SAequiv}.

\begin{cor}
\label{cor:simpS} Suppose $A$ is a coin sequence containing at least one $1$.
Let $A'$ be the coin sequence formed by removing every group of consecutive
$0$'s of even length and concatenating groups of $1$'s when possible. Suppose
$(h_0,h_1,\dots,h_n)$ is the number of consecutive $1$'s in each group of $1$'s
in $A'$. Then $S(A) = h_0 + \sum_{i=1}^{n}(-1)^i h_i - (n \pmod 2)$.
\end{cor}

\begin{prop}
\label{prop:remcomp} A coin sequence $A$ containing at least one heads-up coin
is removable if and only if $S(A)$ is congruent to either $0$ or $1$ modulo
$3$.
\end{prop}

\begin{proof}
We will use induction on the number of coins $k$ in the sequence $A$.  It is
easy to see that the claim holds for $k =1$.  For $k = 2$, we can easily check
the values of $S(A)$ for the possible coin sequences,
\begin{equation*}
S(11) \equiv 2\pmod{3} \qquad\qquad S(10)\equiv 0\pmod{3} \qquad\qquad S(01)
\equiv 1 \pmod{3},
\end{equation*}
\noindent and we see that the claim holds for $k = 2$.

Assume that the claim holds for some $k \geq 2$.  Let $A$ be a coin sequence
with length $k + 1$ such that $A$ contains at least one $1$.  If $A$ is
removable then removing some coin in $A$ yields another coin sequence $A'$ with
length $k$ that is also removable (and so contains at least one $1$). By the
induction hypothesis, $S(A')$ is congruent to $0$ or $1$ modulo $3$ and so by
Lemma~\ref{lem:SAequiv}, $S(A)$ is congruent to $0$ or $1$ modulo $3$ as well.

Now suppose $A$ is not removable.  Then for any choice of coin to remove, the
resulting coin sequence $A'$ is again not removable.  For some coin in $A$,
removing it will leave at least one heads-up coin remaining in $A'$ and so, by
the induction hypothesis, $S(A') \equiv 2 \pmod{3}$.  Applying
Lemma~\ref{lem:SAequiv} again we have that $S(A) \equiv 2 \pmod{3}$.  Thus the
claim holds for coin sequences of length $k+1$ and so by induction the
proposition is true.
\end{proof}

Table~\ref{tab:examples} shows some example applications of
Corollary~\ref{cor:simpS} and Proposition~\ref{prop:remcomp}.

\begin{table}[h]
\begin{center}
\begin{tabular}{lllcc}
    Sequence        & \qquad\qquad $S(A)$\hfill &       & Removable?\\
    \hline
    100110         & $3-0-1$&$\equiv 2 \pmod{3}$     & No\\
    000101         & $0-1+1$&$\equiv 0 \pmod{3}$     & Yes\\
    000110011001   &  $-5$ &$\equiv 1 \pmod{3}$     & Yes\\
    001001100101    &  $4-1-1$&$\equiv 2 \pmod{3}$   & No\\
    111001100001110010011    & $11$&$\equiv 2\pmod{3}$  & No\\
    000110101110100110110    & $-2+1-3+3-2$ &$\equiv 0\pmod{3}$     & Yes
\end{tabular}
\caption{Applications of Corollary~\ref{cor:simpS} and
Proposition~\ref{prop:remcomp}} \label{tab:examples}
\end{center}
\end{table}

\begin{rem} The parity sum provides an effective (polynomial time) means of determining
the removability of sequences.  The proof of Proposition~\ref{prop:remcomp}
shows that there is also an effective algorithm for removing the coins of a
removable sequence: remove {\em any} heads-up coin that either removes all of
the coins or leaves at least one coin heads-up.  Since at least one of the
first two heads-up coins from the left end of the sequence can be removed in
this way, we can remove the coins in polynomial time.
\end{rem}

\section{The No Gaps Variation and a Finite Automaton}

Let $R$ denote the set of all sequences that are removable in the linear
arrangement without gaps variation of the coin-removal problem.  We saw in
Proposition~\ref{prop:remcomp} that we can test a sequence $A$ for membership
in $R$ by calculating the parity sum $S(A)$.  We can also test $A$ for
membership in $R$ by using a {\em deterministic finite automaton} (DFA), which
we present in this section.  Using the DFA, we can then determine a recursive
formula for the number of sequences in $R$ of length $n$.

\begin{figure}[h]
\begin{center}
\includegraphics{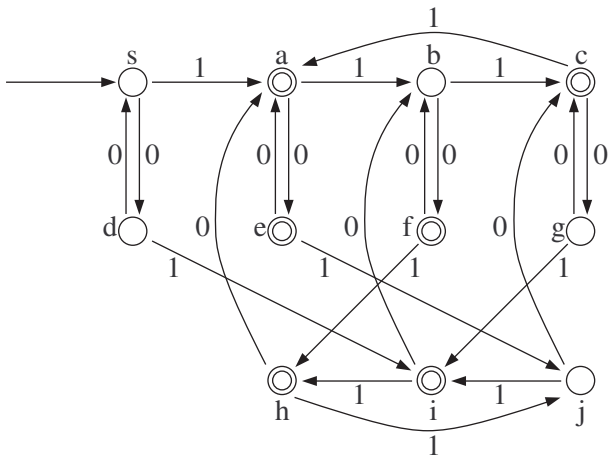}
\caption{DFA to recognize the set $R$ of removable sequences}
\label{fig:dfa_rem}
\end{center}
\end{figure}

Figure~\ref{fig:dfa_rem} depicts a DFA.  Given a sequence $A$ of $0$'s and
$1$'s, we test $A$ for acceptance by the DFA in the following manner:
\begin{itemize}
\item Start at the node labelled $s$, the {\em start} node.
\item Read the symbols of $A$, one at a time, from left to right.  As a symbol
is read, follow the corresponding arrow to the next node (possibly the same
node).
\item Continue reading symbols and following arrows.  After finishing the
sequence, if we are at a node marked with a double circle (an {\em accept}
node) then the DFA accepts the sequence, otherwise it rejects the sequence.
\end{itemize}

\noindent For example, the sequence $1010011$ would take us along the following
path of nodes:
\begin{equation*}
s \overset{1}{\rightarrow} a \overset{0}{\rightarrow} e
\overset{1}{\rightarrow} j \overset{0}{\rightarrow} c \overset{0}{\rightarrow}
g \overset{1}{\rightarrow} i \overset{1}{\rightarrow} h.
\end{equation*}
Node $h$ is an accept node so our DFA will accept the sequence $1010011$.

In fact, the DFA in Figure~\ref{fig:dfa_rem} will {\em recognize}
$R$---sequences will be accepted by the DFA if and only if they are elements of
$R$. To see this, recall that $S(A)$ is computed by adding or subtracting $1$'s
depending on the (modified) parity of the groups of $0$'s separating the groups
of $1$'s. We then compute $S(A) \pmod{3}$ to determine if $A$ is in $R$ or not.
Our DFA performs these computations for us. The nodes $a$, $b$ and $c$ are used
when we are adding $1$'s; they correspond to values of $S(A)$ congruent to $1$,
$2$ and $0$ modulo $3$, respectively. We accept when $S(A) \equiv 0$ or $1$
modulo $3$ (nodes $a$ and $c$).  Nodes $h$, $i$ and $j$ are used when we need
to subtract $1$'s and also keep track of the value of $S(A)$ modulo $3$.  Since
in this case we need to subtract an extra $1$ to complete the computation of
$S(A)$, we only accept at nodes $h$ and $i$. Nodes $d$, $e$, $f$ and $g$ are
used to change from addition to subtraction.

A set of sequences that is recognized by a DFA is called a {\em regular}
language.  Regular languages form the simplest kind of languages studied in
formal language theory.  Other languages, such as palindromes (sequences that
read the same forward and backward) and those containing sequences of the form
$1^p$ for prime $p$ are too complex to be recognized by a DFA.  For these we
need more powerful machines such as {\em pushdown automata} and {\em Turing
machines}.  More information on languages and automata theory can be found in
\cite{hopcroftullman} and \cite{linz}.

\begin{rem}
The results above are analogous to results obtained for Peg Solitaire.  It is
known that the set of configurations that are solvable, i.e., can be reduced to
a single peg, also forms a regular language \cite{mooreeppstein,ravikumar}.
\end{rem}

\begin{figure}[h]
\begin{center}
\includegraphics{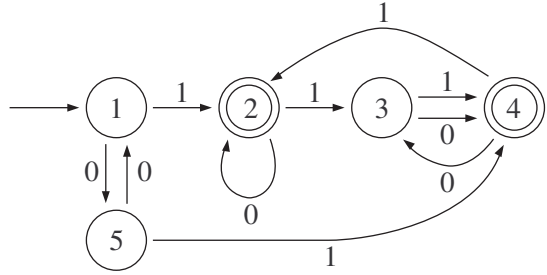}
\caption{Minimized DFA to recognize the set $R$ of removable sequences}
\label{fig:dfamin_rem}
\end{center}
\end{figure}

We can use a state-minimization algorithm \cite[p. 62]{linz} to eliminate
redundant states in our DFA.  The equivalent minimized DFA is shown in
Figure~\ref{fig:dfamin_rem} (the start node is node $1$).  We will use this
presentation to obtain a formula for the number of sequences in $R$ of length
$n$, which corresponds to the number of paths of length $n$ from node $1$ to
either node $2$ or node $4$ (the accept nodes).

Let $M$ be the matrix $\left(M_{ij}\right)$ where entry $M_{ij}$ is the number
of edges from node $i$ to node $j$ in Figure~\ref{fig:dfamin_rem}:
\begin{equation*}
M =
\begin{pmatrix}
    0 & 1 & 0 & 0 & 1\\
    0 & 1 & 1 & 0 & 0\\
    0 & 0 & 0 & 2 & 0\\
    0 & 1 & 1 & 0 & 0\\
    1 & 0 & 0 & 1 & 0\\
\end{pmatrix}.
\end{equation*}

\noindent The matrix $M$ is called the {\em adjacency matrix} of the DFA.  The
adjacency matrix has a very nice property: for any $n \geq 1$, the $(i,j)$
entry in $M^n$ is the number of paths of length $n$ from node $i$ to node $j$.
Thus the number of sequences in $R$ of length $n$ is simply the sum of entries
$(1,2)$ and $(1,4)$ in $M^n$.

Let $r_n$ denote the number of sequences of length $n$ in $R$. Using a computer
to calculate $M^n$ from $n=1$ to $n=7$, we find that the first several terms of
$r_n$ are $1,2,5,10,21,42$ and $85$. These terms appear in sequence A000975 of
the On-line Encyclopedia of Integer Sequences \cite{integerseq}.  The
Encyclopedia entry for A000975 gives a recursive formula for the sequence that
we now show describes $r_n$.

\begin{prop}
The sequence $\{r_n\}$ satisfies the recurrence relation $r_{n+1} = r_n +
2r_{n-1} + 1$.
\end{prop}
\begin{proof}
Let the first row of $M^n$ be denoted by $(a_n,b_n,c_n,d_n,e_n)$.  Then $r_n =
b_n + d_n$.  From matrix multiplication we see that $a_n$, $b_n$, $c_n$, $d_n$
and $e_n$ satisfy the recurrence relations
\begin{alignat*}{6}
    a_{n+1} &= e_n               &\qquad\qquad c_{n+1} &= b_n + d_n &\qquad\qquad & e_{n+1} &= a_n \\
    b_{n+1} &= a_n + b_n + d_n   &\qquad\qquad d_{n+1} &= 2c_n + e_n
\end{alignat*}
with $a_1 = c_1 = d_1 = 0$ and $b_1 = e_1 = 1$.  We will show that $b_{n+1} =
b_n + 2b_{n-1}$ and $d_{n+1} = d_n + 2d_{n-1} + 1$. Notice that $a_n = e_{n-1}
= a_{n-2}$ and $a_n + a_{n-1} = 1$.  Then
\begin{align*}
2b_{n-1} &= 2a_{n-2} + 2(b_{n-2} + d_{n-2})\\
        &= a_n + e_{n-1} + 2c_{n-1}\\
        &= a_n + d_n\\
\end{align*}
and so $b_n + 2b_{n-1} = a_n + b_n + d_n = b_{n+1}$.  Further, the above shows
that $d_n = 2b_{n-1} - a_n$.  Thus
\begin{align*}
d_n + 2d_{n-1} + 1  &= d_n + 4b_{n-2} - 2a_{n-1} + 1\\
                    &= d_n + 2b_n - 2b_{n-1} - 2a_{n-1} + 1\\
                    &= 2b_{n-1} - a_n + 2b_n - 2b_{n-1} - 2a_{n-1} + 1\\
                    &= 2b_n - a_{n-1}\\
                    &= 2b_n - a_{n+1}\\
                    &= d_{n+1}.
\end{align*}
Finally, from the fact that $r_n = b_n + d_n$, we have immediately that
$r_{n+1} = r_n + 2r_n + 1$ and our conjectured formula is proved correct.
\end{proof}

\begin{rem}
The sequence $\{d_n\}$ represents the number of sequences that end at accept
state $4$ and satisfies $d_{n+1} = r_n$ for $n \geq 1$. The sequence $\{b_n\}$
represents the number of sequences that end at accept node $2$; it is the
Jacobsthal sequence (A001045 of \cite{integerseq}).
\end{rem}

\section{Circular Arrangements}
Suppose now that we arrange the coins in a circle.  Again we will only remove
heads-up coins and flip any adjacent coins.  As with the linear arrangements,
we can consider two ways of handling the empty space created by removing a
coin: leave the gap or push the coins together to remove any gaps.

First consider the circular arrangement with gaps variation. (This is actually
a {\em Monthly} problem number 10459, ``Solitaire on a Circle'',
\cite{monthlyprob}).  An example of a removable arrangement in this variation
is $10010$:

\begin{center}
\includegraphics{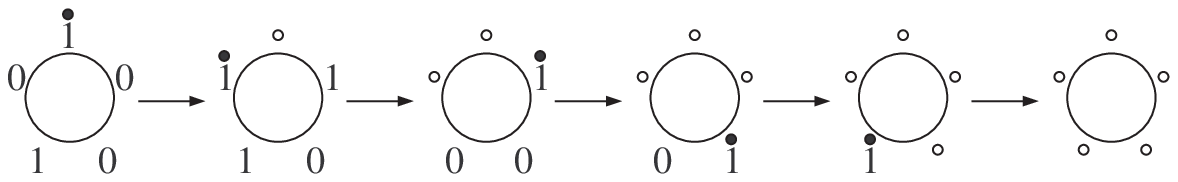}
\end{center}

\noindent We can easily see that the only removable arrangements of lengths $1$
or $2$ are $1$, $10$ and $01$.  For $3$ or more coins, the first coin we remove
will create a {\em linear} arrangement (with gaps).  Also, since the coins are
arranged in a circle, the coin we remove must be the middle coin of some group
of three: $010$, $011$, $110$ or $111$.  In any case, removing the middle coin
changes the number of heads-up coins from even to odd or odd to even.   By
applying Proposition~\ref{prop:linegaps}, we immediately have the following
result.

\begin{prop}
For $n\geq 3$, a circular arrangement of $n$ coins with gaps is removable if
and only if the number of heads-up coins is positive and even.
\end{prop}

Now suppose that when a coin is removed, the remaining coins are pushed
together to eliminate any gaps.  For example, the arrangement $10010$ is
removable, as the following sequence of moves shows:

\begin{center}
\includegraphics{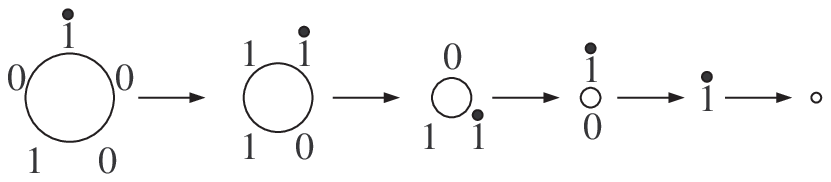}
\end{center}

The only interesting arrangements in this case contain at least two coins.  For
this variation, the parity of the tails-up coins provides us with an
appropriate invariant.

\begin{prop}
For $n\geq 2$, a circular arrangement of $n$ coins without gaps is removable if
and only if it contains at least one heads-up coin and the number of tails-up
coins is odd.
\end{prop}
\begin{proof}
We use induction on the number $k$ of coins.  Clearly for $k=2$ the proposition
is true---the only removable arrangements are $01$ and $10$ (which are actually
the same arrangement).

Assume that the proposition is true for $k$ coins with $k \geq 2$.  Let $A$ be
an arrangement of $k+1$ coins containing at least one heads-up coin.  Any coin
that is removed to obtain $A'$ must be the middle coin of a group of three,
either $010$, $011$ $110$ of $111$.  In each case, removing the middle coin to
obtain an arrangement of $k$ coins will not change the parity of the number of
tails-up coins.

Thus if $A$ has an even number of tails-up coins then $A'$ does as well and by
the inductive hypothesis, $A'$ is not removable.  This holds for any choice of
coin to remove from $A$, so $A$ is not removable (otherwise at least one $A'$
would itself be removable). If the number of tails-up coins in $A$ is odd then
$A'$ also has an odd number of tails-up coins and is removable.  Therefore $A$
is also removable and our proposition is proved.
\end{proof}

\section{Future Investigations---Other Arrangements and a Game}

Other ways to arrange the coins naturally suggest themselves.  For example, we
might place the coins in a rectangular arrangement and flip the horizontal and
vertical neighbors of a removed coin.  An example of a removable rectangular
configuration (for a $2$ by $4$ rectangle) is
\begin{center}
\begin{tabular}{cccc}
    1 & 0 & 1 & 0\\
    0 & 1 & 0 & 1
\end{tabular}
\end{center}
while
\begin{center}
\begin{tabular}{cccc}
    0 & 1 & 1 & 0\\
    0 & 0 & 0 & 0
\end{tabular}
\end{center}
is an example of a non-removable configuration.  Even more generally, we might
place the coins at the vertices of a (connected) graph and on removing a coin,
flip any coins at adjacent vertices.  These situations are similar to those
considered in ``Lights Out'' games.  In a ``Lights Out'' game, lights are
arranged in a rectangular lattice---turning off or on one light will toggle the
state of the adjacent lights.  However, the light is {\em not} removed from the
lattice, unlike the coin removal games considered here.  For more information,
see \cite{andersonfeil}.

Siler \cite{siler} has examined the case of rectangular grids and showed that
for any $m \times n$ grid, if $n$ is odd then $G$ is removable if and only if
$G$ contains an odd number of heads-up coins.  Further, he is able to
completely characterize the case of rectangular grids of size $2 \times n$:

\begin{thm}[Siler] Let $G$ be a $2 \times n$ grid of coins.  If $n = 2k$ is even
then $G$ is removable if and only if the number of heads-up coins in $G$ is
even and nonzero and $G$ is not of the form
\begin{align*}
&0 a_1 a_3 \cdots a_{2k - 5}0\\
&0 a_2 a_4 \cdots a_{2k - 4}0
\end{align*}
where $a_i \in \{00, 11\}$.  If $n$ is odd then $G$ is removable if and only if
$G$ contains an odd number of heads-up coins.
\end{thm}

It would also be interesting to investigate a two-player game based on the
original coin-removal problem (linear arrangement with gaps).  For the
two-player game, players take turns removing coins according to the usual
rules; the loser is the player who cannot pick up a coin. For example, with the
initial arrangement of $111$, the first player can win by selecting the middle
coin, which leaves $0\circ 0$ so the second player has no allowed move. But
with the initial arrangement of $101$, the second player can win no matter what
the first player does: after the first player's move, the arrangement will be
$\circ 11$ which the second player can reduce to $\circ \circ 0$.  Of
particular interest are the formal language properties of the set of winning
sequences for one player or the other.  A similar game based on Peg Solitaire
was suggested in \cite{ravikumar} and examined in \cite{mooreeppstein}.

\end{document}